# Arnold's Weak Resonance Equation as the Model of Greek Ornamental Design

**Faina Berezovskaya[1], Georgiy Karev[2]**

[1] **Department of Mathematics, Howard University, Washington, DC 20059, USA**

Email: fberezovskaya@howard.edu

[2] **National Centre for Biotechnology Information, NIH, Bethesda, MD 20894, USA**

Email: karev@ncbi.nlm.nih.gov

**Abstract**

We propose and study a mathematical model that qualitatively reproduces several ancient ornamental designs than one can see in historical museums of Crete and Athens. The designs contain several rings that circumscribe a fixed number of "flowers" (centers or spirals), specific to each design.

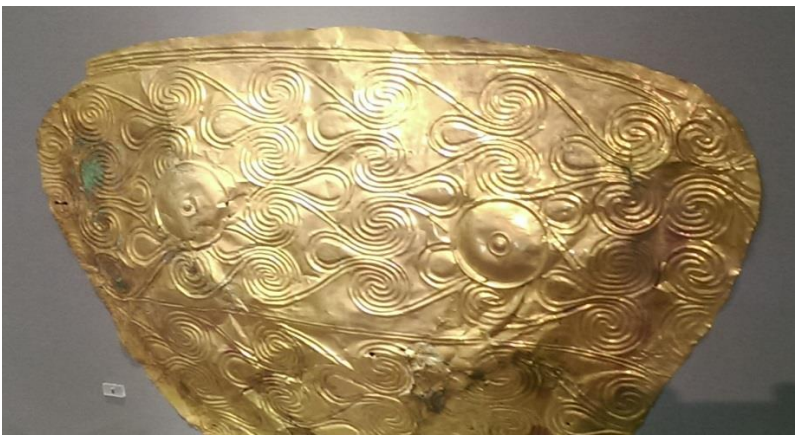

The model is based on a complex differential equation of "weak resonance" (Arnold 1977). We analyze the role of the model parameters in giving rise to different peculiarities of the repeated designs, in particular, the "dynamical indeterminacy". The model allows tracing design changes under parameter variation, as well as to construct some new ornamental designs. We discuss how observed ornamental design may reflect some philosophical ideas of ancient inhabitants of Greece.

**Introduction**

*1.1. The model as an equation with complex variable*

Clear patterns exist in ancient ornamental designs, such as the ones that one can see in historical museums of Crete and Athens (see examples in Appendix 3). The designs contain bands of different but fixed numbers of "flowers" (mathematically, centroids or spirals), "spider nets" and various types of cycles up to "stars" (see below the definitions). Boundaries of the bands may be smooth or star-shaped curves; the bands are "connected" by smooth lines. These designs reminded us of some phase portraits of quasi-Hamiltonian dynamical systems, invariant under rotations by the specific angles $2\pi/n$ with integer $n$ (see Appendix 3, Figures 2,3,5).

The aim of this paper is to find the simplest mathematical model that can describe key features of these ornaments. To this end, we consider the equation proposed by V. Arnold for analysis the problem "Loss of stability of self-sustained oscillations" (see [2], [3] and references therein).

The Arnold' equation is a complex differential equation that describes an equivariant vector field, i.e. symmetric vector field invariant with respect to rotation to the angle $2\pi/n$

$$z' = \epsilon z + zA(|z|^2) + B\bar{z}^{n-1} \qquad (1.1)$$

where $z = x + iy$ is a point in the complex plane, integer $n \geq 4$, $\epsilon = \epsilon_1 + i\epsilon_2$; the function $A(|z|)$ is given by the formula

$$A(|z|) = \sum_{k=1}^{s} A_k |z|^{2k}, \quad A_k = A_k^1 + iA_k^2, B = B_1 + iB_2, \qquad (1.2)$$

where integer $s = \frac{n}{2} - 1 = \frac{n-2}{2}$ if $n$ is even, $s = \frac{n-1}{2} - 1 = \frac{n-3}{2}$ if $n$ is odd. Evidently, equation (1.1) has equilibrium $z = 0$ for any $A_k$, $B$, and $\epsilon_1, \epsilon_2$ are its eigenvalues; the equation can have also "peripheral" equilibria $E(z^*)$ if for some $z^* \neq 0$, $\epsilon z^* + z^* A(|z^*|^2) + B\bar{z^*}^{n-1} = 0$.

Originally, Equation (1.1) was constructed for describing the loss of stability of self-oscillations in maps; Equation (1.1) is an approximation of maps under *"strong resonance"* if *n*=1,2,3,4 and *"weak resonance"* if *n*≥ 5. Cases for $n \leq 4$ were studied in [2-9] (and discussed in many other works (see for example [10]) using both analytical and numerical methods of bifurcation theory. To the best of our knowledge, cases of $n \geq 5$ are still not investigated completely.

The main problem considered in [2,3] was description of phase –parameter portraits of Equation (1.1) in a neighborhood of co-dimension 2 bifurcation for $\epsilon = 0$. Specifically, the goal was to reveal sequences of all co-dimension 1 bifurcations that arise in the vicinity of $(0, \epsilon_2)$, for any fixed coefficients $A \neq 0, B = 1$. Notably, all cases $n \leq 4$ demonstrate essentially different phase behaviors, and corresponding equations have different "organizing centers".

In what follows we consider cases for $\underline{n \geq 4}$ and show that Equation (1.1) demonstrates *different kinds* of phase behaviors depending on whether *n* is even or odd. We analyze the role of parameters *B* and $A_s$ for even or odd *s* in the genesis of patterns and repeated designs for different *n*.

Phase portraits of the equation with $n \geq 4$ have patterns that mimic the qualitative features of some of Greek ornamental designs. Figures of Appendix 2 show several phase portraits of Equation (1.1) for different *n*; all of these portraits contain some of the four main types of annular patterns observed in Greek ornamental designs.

***Definitions*** *(see Figure 1.1)*

- *"**n- cycle**" is a separatrix limit cycle composed of n saddles "connected" by their separatrices; we call it an "n-star" if it is not convex and "convex n-cycle" otherwise;*
- *"**centroid**" is a pattern composed of a center or spiral equilibrium together with a set of orbits around it;*
- *"**n-flower ring**" is a pattern consisting of n centroids and n saddles together with their separatrices; each centroid is placed inside the "leaf" composed by a separatrix cycle or separatrix loop;*
- *"**spider-net**" is a pattern consisting of one or two neighboring saddles together with their separatrices (outgoing to infinity) and hyperbolic-shape orbits between the separatrices.*

Next, we show below that for a wide range of parameters, Equation (1.1) is Hamiltonian, and its phase portraits represent certain collections (unions) of patterns described above. Detailed descriptions of phase portraits of Hamiltonian system are given in s.3.

We also present some portraits, which cannot be reproduced by Hamiltonian equations but can be considered as portraits of Hamiltonian equation under small variations of its coefficients.

**Statement A.** *For any integer $n \geq 4$ and small values of parameters $\{\epsilon_2 \geq 0, A_k^1, k = 1, ..., s\}$, the phase-parameter portrait of equation (1.1) can contain the four patterns defined above.*

In Discussion, we try to explain and comment on the "dynamical indeterminacy" that to us is the most interesting peculiarity of the Greek designs. In terms of dynamical systems, it implies the presence of domains in the phase plane such that the system shows essentially different limiting behaviors for close initial values. For example, we observe such domains in Figures A2.7, 8, where trajectories from the neighborhood of the origin (**b**) or the limit cycle (**a**) can reach any of the peripheral points.

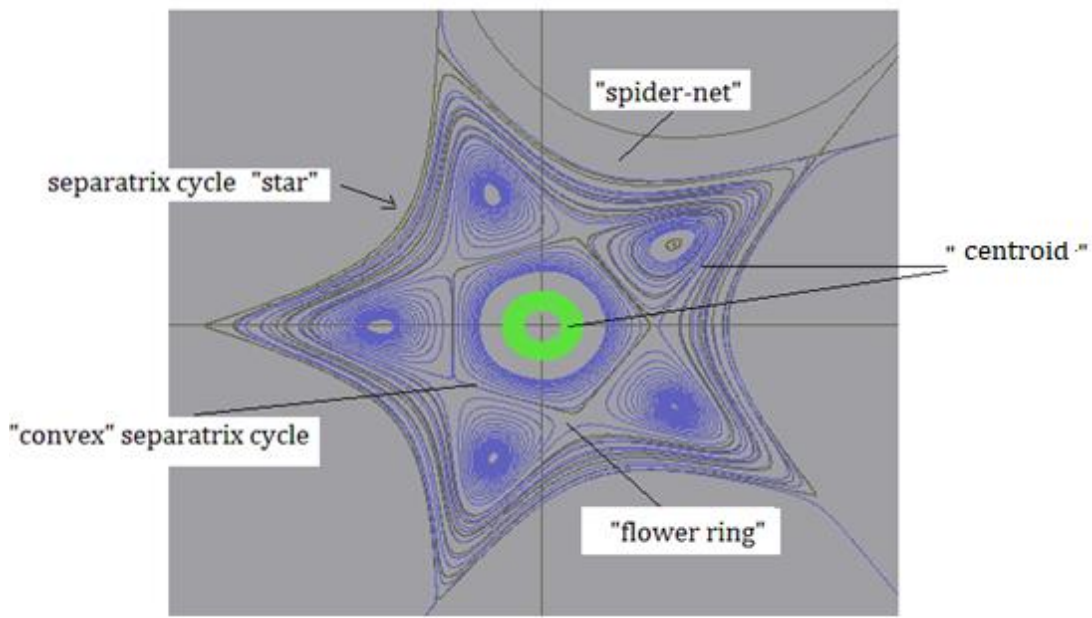

**Figure 1.1.** Example of model (1.1) with *n*=5. The phase portrait contains centroids around *O* and five peripheral equilibria; the peripheral equilibria compose one "5-flower ring" bounded by 2 separatrix cycles; any two neighboring saddles contain "spider-nets".

2. **Some properties of the model**

   *2.1. The equation in different coordinate systems*

In order to analyze Equation (1.1), it is convenient to cast it in both polar and the Descartes coordinates.

In polar coordinates $(r, \varphi)$: $z = re^{i\varphi}$ ($x = r\cos(\varphi), y = r\sin(\varphi)$), Equation (1.1) becomes the following System:

$$r' = r(\epsilon_1 + \sum_{k=1}^{S} A_k^1 r^{2k} + r^{n-2}(B_1 \cos(n\varphi) + B_2 \sin(n\varphi)), \qquad (2.1)$$

$$\varphi' = \epsilon_2 + \sum_{k=1}^{S} A_k^2 r^{2k} + r^{n-2}(-B_1 \sin(n\varphi) + B_2 \cos(n\varphi)).$$

In $(x, y)$ coordinates Equation (1.1) becomes the following System:

$$x' = \epsilon_1 x - \epsilon_2 y + \sum_{k=1}^{S}(A_k^1 x - A_k^2 y)(x^2 + y^2)^k + (B_1 p(x, y) + B_2 q(x, y)), \qquad (2.2)$$

$$y' = \epsilon_2 x + \epsilon_1 y + \sum_{k=1}^{S}(A_k^2 x + A_k^1 y)(x^2 + y^2)^k + (-B_1 q(x, y) + B_2 p(x, y)),$$

where $p(x, y)$, $q(x, y)$ are polynomials of $(n-1)$-th order, such that $p_x(x, y) = q_y(x, y), p_y(x, y) = q_x(x, y)$.

Rotating the phase plane by angle $\beta$ such that $\cos\beta = \dfrac{B_1}{\sqrt{B_1^2 + B_2^2}}$, $\sin\beta = \dfrac{B_2}{\sqrt{B_1^2 + B_2^2}}$, we can rewrite Systems (2.1) and (2.2), correspondingly, as follows:

$$r' = r(\epsilon_1 + \sum_{k=1}^{S} A_k^1 r^{2k} + B\, r^{n-2} \cos(n\varphi)) \equiv P(r, \varphi) \qquad (2.3)$$

$$\varphi' = \epsilon_2 + \sum_{k=1}^{S} A_k^2 r^{2k} - B r^{n-2} \sin(n\varphi) \equiv Q(r, \varphi)$$

where $B = \sqrt{B_1^2 + B_2^2}$, and

$$x' = \epsilon_1 x - \epsilon_2 y + \sum_{k=1}^{S}(A_k^1 x - A_k^2 y)(x^2 + y^2)^k + B\, p(x, y) \equiv F_1(x, y), \qquad (2.4)$$

$$y' = \epsilon_2 x + \epsilon_1 y + \sum_{k=1}^{S}(A_k^2 x + A_k^1 y)(x^2 + y^2)^k - B q(x, y) \equiv F_2(x, y).$$

Notice that

$$P(r, \varphi) = F_1(x, y) \cos(\varphi) + F_2(x, y) \sin(\varphi), \qquad (2.5)$$

$$Q(r, \varphi) = (F_2(x, y) \cos(\varphi) - F_1(x, y) \sin(\varphi))/r$$

where $x = r\cos(\varphi), y = r\sin(\varphi)$.

Equation (1.1), as well Systems (2.3) and (2.4) possess many specific properties, some of which we discuss below.

*2.2. Case B=0.*

If *B*=0, then Equation (1.1) reads

$$z' = \epsilon z + zA(|z|^2),  \quad (2.6)$$

where the function *A* is defined by (1.2).

This differential equation serves as an approximation of the Poincaré map, which was applied to analysis of loss of stability of a closed orbit (limit cycle) [1, 10]. Let us recall some important properties of Equation (2.6).

In polar coordinates, Equation (2.6) becomes

$$r' = r(\epsilon_1 + \sum_{k=1}^{s} A_k^1 r^{2k}) \equiv r\, P_1(r), \quad (2.7)$$

$$\varphi' = \epsilon_2 + \sum_{k=1}^{s} A_k^2 r^{2k} \equiv Q_1(r).$$

In particular, for small $\epsilon_1$ and $\epsilon_2 \neq 0, A_1^1 \neq 0$; truncation of the system:

$$r' = r(\epsilon_1 + A_1^1 r^2), \varphi' = \epsilon_2$$

serves as a model system for the Andronov-Hopf bifurcation of changing stability of equilibrium O(0,0), where $A_1^1$ is the first Lyapunov value. The bifurcation is supercritical, accompanied by appearance/disappearance of a *stable* limit cycle if $\{A_1^1 < 0, \epsilon_1 > 0\}$, and subcritical if $\{A_1^1 > 0, \epsilon_1 < 0\}$, accompanied by appearance/ disappearance of an *unstable* limit cycle. If $A_1^1 = 0$, then the number and stability of limit cycles is defined by the second Lyapunov value $A_2^1$; if $A_2^1 \neq 0, \epsilon_2 \neq 0$, then the next truncation system, containing term $A_2^1 r^4$, describes the generalized Andronov-Hopf bifurcation.

Thus, the first Equation (2.7) defines the number and stability of limit cycles of System (2.7).

For $\epsilon_1 \cong 0$, equilibrium O is a "weak spiral" (a center in linear approximation); the direction of orbit rotation of (2.7) close to $O$ is defined by the sign of $\epsilon_2$, and in a general case, by the sign of $Q_1$. The direction of rotation can change for $r = r^{**}$, where $r^{**}$ is a root of the function $Q_1$, so that $\varphi' = 0$.

The point $(r^*, \varphi)$ such that $P_1(r^*) \neq 0, Q_1(r^*) = 0$ is called a "*quasi-equilibrium*" of (2.7). Quasi-equilibria are always composed of a circle with radius $r^*$ such that $\varphi' = Q_1(r^*) = 0$; we call to this circle a *quasi-equilibrium cycle*. In the vicinity of every quasi-equilibrium the orbits of the system change direction of rotation (see Figure 2.1). Taking into consideration that limit cycles of (2.7) correspond to the roots of polynomial $P_1$, and quasi-equilibrium cycles correspond to the roots of polynomial $Q_1$ and applying the Descartes' Rule of Signs ("The number of positive roots of the polynomial is either equal to the number of sign differences between consecutive nonzero coefficients, or is less than it by an even number") we prove the following statement:

**Proposition 1.** *The number $n_c$ of limit cycles and the number $n_q$ of quasi-equilibrium cycles of System (2.7) does not exceed s, where s is defined in (1.2).*

.

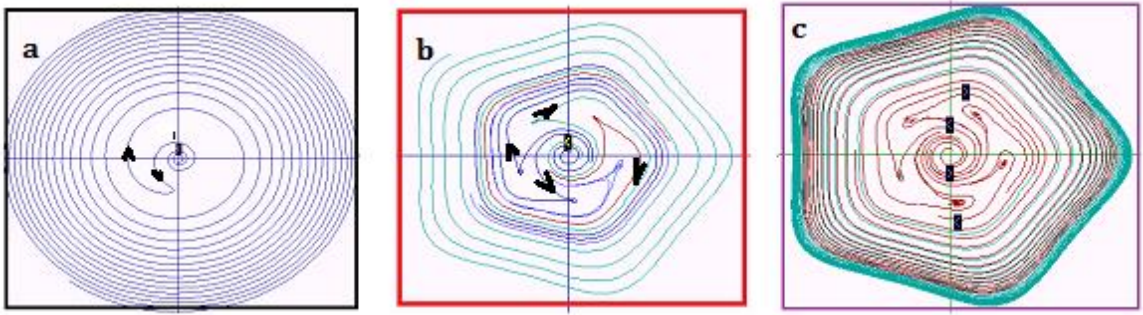

Figure.2.1. Quasi-equilibria of System (2.9), with $n=5$, $\varepsilon_1 = .0001, \varepsilon_2 = .1, A_1^2 = -.1$. Equilibrium $O$ is "weakly unstable" in all three panels. **a**: quasi-equilibrium point as $A_1^1 = 0, B = 0$; **b**: appearance of peripheral equilibria for $A_1^1 = 0, B = .01$; **c:** appearance of a limit cycle and peripheral equilibria for $A_1^1 = -.01, B = .01$.

Let us study the role of the term $B\bar{z}^{n-1}$ in equation (1.1). In Descartes coordinates it becomes $Bp(x, y)$ in the 1st equation of System (2.4), and $-Bq(x, y)$ in the 2nd equation of (2.4). Here, $p(x, y)$, $q(x, y)$ are polynomials of ($n$-1)-th power, and $p_x(x, y) = q_y(x, y), p_y(x, y) = -q_x(x, y)$. In polar coordinates, this term corresponds to $Br^{n-1}\cos(n\varphi)$ in the first equation of System (2.3) and to $-Br^{n-2}\sin(n\varphi)$ in the second equation of (2.3). We show below that these terms are responsible for appearance of "peripheral" equilibria $E(r_k, \varphi_k)$, where $\varphi_k = \pm\frac{\pi}{2n} + 2\pi k/n$, $k$=0,…$n$-1.

3. Hamiltonian model

3.1. Hamiltonian

It is known [1] that a system in the Descartes' coordinates is Hamiltonian if its divergence vanishes. The divergence of Equation (1.1) taken in form (2.4)

$$Div = \frac{\partial F_1}{\partial x} + \frac{\partial F_2}{\partial y} =$$

$$2(\epsilon_1 + (x^2 + y^2)(2A_1^1 + (x^2 + y^2)(3A_2^1 + (x^2 + y^2)(4A_3^1 + (x^2 + y^2)(\ldots)))) =$$

$$2(\epsilon_1 + r^2(2A_1^1 + r^2(3A_2^1 + r^2(4A_3^1 + r^2(5A_4^1 \ldots))))$$

vanishes as

$$\epsilon_1 = 0, A_k^1 = 0, k = 1, \ldots, s. \tag{3.1}$$

Then the following statement is true:

***Proposition 2***. *If conditions* (3.1) *hold, then Equation (1.1) is Hamiltonian and can be written in polar coordinates in the form*

$$r^` = B\,r^{n-1}\cos(n\varphi) = -r\frac{\partial H}{\partial \varphi} \equiv P_H(r, \varphi) \tag{3H}$$

$$\varphi^` = \epsilon_2 + \sum_{k=1}^{s} A_k^2 r^{2k} - Br^{n-2}\sin(n\varphi) = \frac{1}{r}\frac{\partial H}{\partial r} \equiv Q_H(r, \varphi)$$

*with Hamiltonian*

$$H(r,\varphi) = \epsilon_2 \frac{r^2}{2} + \sum_{k=1}^{S} A_k^2 \frac{r^{2k+2}}{2k+2} - B\frac{r^n}{n}\sin(n\varphi) + h, \qquad (3.2)$$

*where h is an arbitrary constant.*

In (*x, y*)- coordinates, Equation (1.1) as given in the form of System (2.4), becomes

$$x` = -(\epsilon_2 y + y \sum_{k=1}^{S} A_k^2 (x^2+y^2)^k + B\, p(x,y) = -\frac{\partial H}{\partial y} \equiv F_{1H}(x,y), \qquad (4H)$$

$$y` = \epsilon_1 y + \sum_{k=1}^{S} A_k^2 y(x^2+y^2)^k - Bq(x,y) = \frac{\partial H}{\partial x} \equiv F_{2H}(x,y)$$

with Hamiltonian

$$H(x,y) = \epsilon_2 \frac{(x^2+y^2)}{2} + \sum_{k=1}^{S} \frac{A_k^2(x^2+y^2)^k}{2k} + BR(x,y) + h, \qquad (3.3)$$

where $R(x,y) = \int \left(q(x,y) + \frac{\partial}{\partial x}\int p(x,y)\,dy\right)dx - \int p(x,y)\,dy$, and *h* is an arbitrary constant.

System (3H) has equilibrium $O(0,0)$ and can have *"peripheral"* equilibria $E^*(r^*,\varphi^*), r^* > 0$, whose coordinates $(r^*,\varphi^*)$ satisfy the system

$$\cos(n\varphi) = 0, \qquad \epsilon_2 + \sum_{k=1}^{S} A_k^2 r^{2k} - Br^{n-2}\sin(n\varphi) = 0. \qquad (3.4)$$

The first Equation of (3.4) defines $2n$ rays $\varphi = \varphi^* = \varphi_j^{\pm}$, where

$$\varphi_j^{\pm} = \pm\frac{\pi}{2n} + \frac{2\pi j}{n}, j = 0,\dots,n-1. \qquad (3.5\pm)$$

Then $\sin(n\varphi_j^{\pm}) = \pm 1$, and coordinates $r^*$ of equilibria $E^*$ are the roots of one of the polynomials

$$P^+ \equiv \epsilon_2 + \sum_{k=1}^{S} A_k^2 r^{2k} - Br^{n-2}, \quad P^- \equiv \epsilon_2 + \sum_{k=1}^{S} A_k^2 r^{2k} + Br^{n-2}, B>0 \qquad (3.6\pm)$$

Let us consider polynomials $P^+, P^-$ along the rays $\varphi_j^+, \varphi_j^-$ for the same fixed $j = 0,\dots n-1$. The corresponding equilibria $E^{\pm}(r^{\pm},\varphi_j^{\pm})$ can be only saddles or centers (see Figure 1.1, A2.1-A2.6) because the system is Hamiltonian.

Clearly, if the equilibrium $E^+(r^+, \varphi_{j_0}^+)$ is a saddle (center) for some $j_0$, then all equilibria $E^+(r^+, \varphi_j^+)$ are saddles (centers), $j = 0, \ldots n-1$. Similar assertion is valid for equilibria $E^-(r^-, \varphi_{j_0}^-)$. These properties allow us to omit index $j$ in the notation of peripheral equilibrium and consider only points $E^+(r^+, \varphi^+)$ and $E^-(r^-, \varphi^-)$ because topological type of equilibrium $E$ does not depend on $j$. For example, if $E^+(r^+, \varphi^+)$ $(E^-(r^-, \varphi^-))$ is a center/a saddle, then the phase plane contains $n$ centers/saddles with the same coordinates $r^+(r^-)$ and $\varphi = \varphi^+ + \frac{2\pi j}{n}, j = 0, \ldots, n-1 / \varphi = \varphi^- + \frac{2\pi j}{n}, j = 0, \ldots, n-1$.

The number and characteristics of equilibria essentially depends on whether $n$ is odd or even.

***Proposition 3.*** *Let $n \geq 5$ be odd. Then*

1) *polynomials $P^+, P^-$ have no more than $s = \frac{n-1}{2}$ real positive roots;*
2) *if polynomial $P^+(P^-)$ has M real roots and J of them are positive, then the polynomial $P^-(P^+)$ also has M real roots and M-J of them are positive;*
3) *let $r^*, r_*$ be maximal and minimal roots among all positive roots of polynomials $P^+$ and $P^-$. Then corresponding equilibria $E^*$ and $E_*$ are saddles;*
4) *every two equilibria $E_j^+, E_{j+1}^+$ corresponding to consequent roots of polynomial $P^+$ are saddle/center or center/saddle; same statement is valid for polynomial $P^-$.*

In order to study model behavior as $x, y \to \infty$, it is useful to consider the system on the Poincaré sphere (see for example [1]). The Poincaré sphere is defined by two transformations of Descartes coordinates $(x, y)$, given by formulas $\{ = \frac{1}{x}, w = \frac{y}{x}\}$ and $\{u = \frac{1}{y}, v = \frac{x}{y}\}$.

Analyzing equilibria in the equators of the Poincaré sphere as $u = 0$, we get the following

***Proposition 4.*** *In the equator of the Poincaré sphere, System (3H) with odd n has n equilibria, which are alternating stable and unstable nodes (see Figure 3.1)*

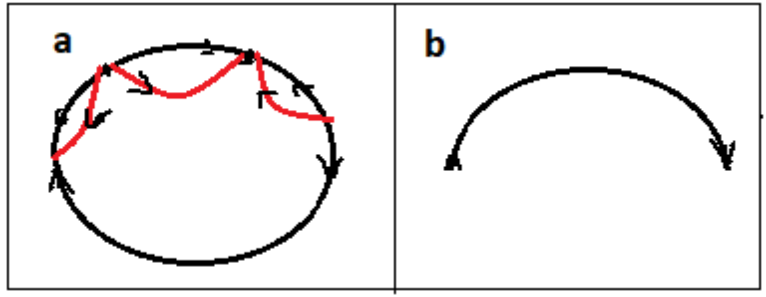

Figure 3.1. Equators of Poincaré sphere, (a) $n$ is odd, (b) $n$ is even, $B < |A_s^2|$.

Using Propositions 3 and 4 we can prove the following

***Theorem 2***. *Let conditions (3.1) hold and assume Equation (1.1) is Hamiltonian. Then for odd $n \geq 5$, the phase portrait of the Equation contains only centroids, n-cycles and spider-nets and can contain no more than $\frac{n-3}{2}$ flower rings.*

**Example 1**. System (3H) with *n*=5 (see Figure 3.2).

In this case, only one ring can exist due to Proposition 3. The right-hand sides of System (3H) are

$$P_H(r,\varphi) = Br^4 \cos(5\varphi), \ Q(r,\varphi) = \epsilon_2 + A_1^2 r^2 - Br^3 \sin(5\varphi), B > 0.$$

Polynomial $P_0 = \epsilon_2 + A_1^2 r^2$ has one positive root if $\epsilon_2 A_1^2 < 0$ and has no real roots if $\epsilon_2 A_1^2 > 0$. In the first case for *B*>0, one of polynomials $P^\pm = \epsilon_2 + A_1^2 r^2 \mp Br^3$ has one real positive root, and second polynomial has one negative and two positive roots. Common number of positive roots of both $P^\pm$ is $\frac{n-1}{2} = 3$, and two of them are saddles. Line $C: 27\epsilon_2 B^2 + 4(A_1^2)^3 = 0$ serves as the boundary between two domains of different phase behaviors.

Let, for example, $P^+$ have two positive roots $r_1^+ = 1, r_2^+ = 2$. Then coefficients of $P^+$ are $\epsilon_2 = -\frac{4B}{3}, A_1^2 = \frac{7B}{3}$. For any $B > 0$ polynomial $P^+$ also has a negative root.

In Figure 3.2 we present the bifurcation diagram of the model for $\epsilon_2 = -4, A_1^2 = 7, B = 3$. Here $P^+/P^-$ has roots $\{r_1{}^+, r_2{}^+, r_3{}^+\} = \{1, 2, -\frac{2}{3}\}, \{r_1{}^-, r_2{}^-, r_3{}^-\} = \{-1, -2, 2/3\}$. Additional examples are given in Appendix 3, Figure A3.2, 3.

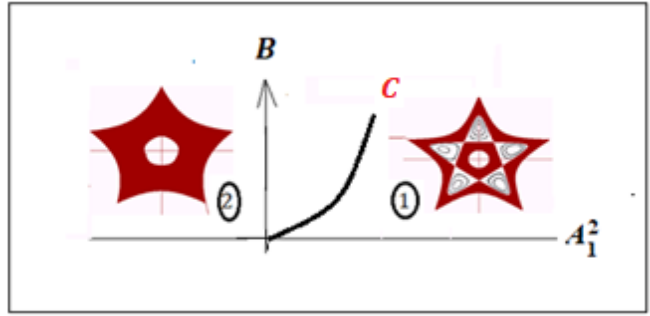

**Figure 3.2.** Bifurcation diagram of Example 1 for n=5, $\epsilon_2 = -4, B = 3$; Domain 1: $A_1^2 = 7$, Domain 2: $A_1^2 = -1$; boundary $C: B = (A_1^2/27)^{3/2}$.

Now let us consider the case, when $n \geq 4$ in Model (3H) is *even*.

Then $s = \frac{n}{2} - 1$ (see (1.2)) and System (3H) can be rewritten in the form

$$r` = B\, r^{n-1} \cos(n\varphi) \tag{3.7}$$

$$\varphi` = \epsilon_2 + \sum_{k=1}^{s-1} A_k^2 r^{2k} + (A_{(n-2)/2}^2 - B\sin(n\varphi))r^{n-2}.$$

Let $E^*(r^*, \varphi^*)$ be a "peripheral" equilibrium of (3.7). Then its $\varphi^*$-coordinate is given by (3.5±), and $r^*$ −coordinate is a positive root of *even* polynomials $P^\pm$,

$$P^\pm \equiv \epsilon_2 + \sum_{k=1}^{s-1} A_k^2 r^{2k} + (A_{\frac{n-2}{2}}^2 \mp B)r^{n-2}. \tag{3.8 ±}$$

Each polynomial $P^\pm$ can have from 0 up to $\frac{n-2}{2}$ positive roots.

Consider even polynomial $P_0 \equiv \epsilon_2 + \sum_{k=1}^{s-1} A_k^2 r^{2k} + A_{\frac{n-2}{2}}^2 r^{n-2}$ that differs from polynomials $P^\pm$ only by the last coefficient. Polynomial $P_0$ has no more than $s = (n-2)/2$ positive roots.

***Proposition 5***. *Let $n \geq 4$ be even, and B>0. Then*

1) polynomials $P^+, P^-$ have no more than $s = \frac{n-2}{2}$ real positive roots;

2) polynomials $P^{\pm}$ have the same number of positive roots as polynomial $P_0$ if $A_{\frac{n-2}{2}}^2 > -B > 0$ and if $0 < B < -A_{\frac{n-2}{2}}^2$;

3) if $B > \left| A_{\frac{n-2}{2}}^2 \right|$, then at least one of polynomials $P^-, P^+$ has a real positive root $r^*$ and this $r^*$ is $r$-coordinate of a saddle equilibrium of (3.7);

4) every two equilibria $E_j^+, E_{j+1}^+$ corresponding to subsequent roots of polynomial $P^+$ are saddle/center or center/saddle; similar statement is valid for the polynomial $P^-$.

Analyzing "infinite" equilibria at the Poincaré sphere (similar to Proposition 4), we get the following

***Proposition 6.*** *In the equator of the Poincaré sphere, System (3.7) has at least 2 equilibria, which are alternating stable and unstable nodes if $B > \left| A_{\frac{n-2}{2}}^2 \right|$ and has no equilibria otherwise (see Figure 3.1b).*

Propositions 4 and 6 are special cases of Proposition 8 given below.

***Theorem 3.*** *Let conditions (3.1) hold and assume Equation (1.1) is Hamiltonian. Then for even $n \geq 4$, the phase portrait of Equation (1.1) contains centroid, n-cycles, and may contain no more than $\frac{n-2}{2}$ flower rings; for $B > \left| A_{\frac{n-2}{2}}^2 \right|$ the phase portrait additionally contains spider nets.*

Proof of the Theorem is given in Appendix 1.

**Example 2**. System (3H) with $n=4$. Here $P_H(r, \varphi) = Br^3 \cos(4\varphi), Q_H(r, \varphi) = 1 + r^2(A_1^2 - B\sin(4\varphi)), B > 0$.

In this case only one ring is possible. Polynomial $P_0 = 1 + A_1^2 r^2$ has no roots when $A_1^2 > 0$, and has one positive root when $A_1^2 < 0$. Polynomials $P^{\pm} = 1 + r^2(A_1^2 \mp B)$ have no roots if $0 < A_1^2 + B$ and $A_1^2 - B > 0$, correspondingly. They can have two positive roots if $A_1^2 < 0, B < -A_1^2$ and only one root if $|B| > A_1^2$ (see Figures 3.3, A2.1 and A2.3).

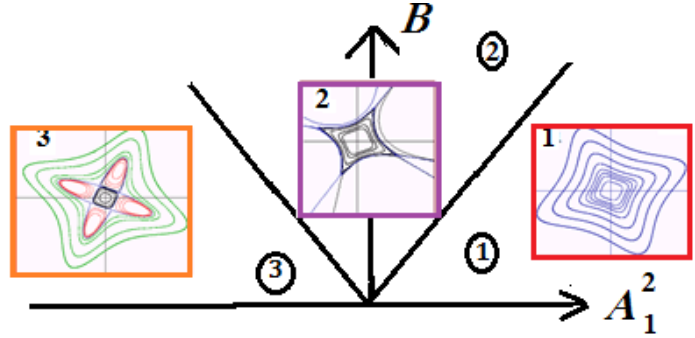

**Figure 3.3.** Bifurcation diagram of the model in Example 2. $n=4$, $\epsilon_1 = 0, \epsilon_2 = 1$. Domain 1: $A_1^2 = 1$, $B = .7$; Domain 2: $A_1^2 = 1, B = 1.2$; Domain 3: $A_1^2 = -1, B = .7$. Boundaries between Domains: $B = \pm A_1^2$.

*3.2. Rearrangements of phase portraits of Hamiltonian model*

Now we consider the how "repeated structures" appear in a phase plane of Hamiltonian Equations (1.1) and (3.1), assuming for simplicity that $\epsilon_2 > 0$. Firstly, let us assume that all coefficients $A_k^2 > 0$. For any $n$, phase curves starting in the vicinity of the equilibrium $O$ result in a "centroid", i.e. a family of closed cycles. Next, we need to differentiate between the cases of odd and even $n$.

Let $n$ be *odd*. Then there exists a "separatrix cycle" for any fixed $B \neq 0$ and some $A_k^2$ that serves as a boundary of a family of closed curves due to Proposition 3 (see Figures 1.1, A3.2, A3.4-A3.6). The equilibria in the equator of Poincaré sphere are stable and unstable nodes (Proposition 4). So, the phase plane contains "spider-nets" (see Figure 1.1). Notice now that the described picture is complete if one of polynomials $P^\pm$ has only one positive root; it happens, for example, if all coefficients $A_k^2$ are positive. If for the same $B$ some of $A_k^2$ are negative, then the "flower-ring" (a ring of centroids inside the separatrix cycles) can appear. Notice, that the number of "flower-rings" cannot be more than $\frac{n-3}{2}$. The appearance of any such ring is accompanied by change of the sign of coefficient $A_k^2$ for some $k$; it implies change of a number of alternating sign in sequence $\{A_k^2, B\}$, $k=1,...s$.

Remark, that structure of phase portraits for odd $n$ depends also on whether $s = \frac{n-3}{2}$ is odd or even. For odd $s$ (see Figures A2.4, A2.6, left panels as n=7, 11), "specific" flower-rings can arise, which are similar to some Greek designs (Figure A2.5).

Let $n$ be *even*. Then the number of flower rings, as well as the number of centers is no more than $\frac{n-2}{2}$. Let $\epsilon_2 > 0$. If $0 \leq B < \left|A^2_{\frac{n-2}{2}}\right|$, then polynomials (3.8±), $P^\pm \equiv \epsilon_2 + \sum_{k=1}^{s-1} A^2_k r^{2k} + (A^2_{\frac{n-2}{2}} \mp B) r^{n-2}$, and $R \equiv \epsilon_2 + \sum_{k=1}^{s-1} A^2_k r^{2k} + A^2_{\frac{n-2}{2}} r^{n-2}$ have the same number of sign changes in their sequences of coefficients. If polynomial $R$ has no real roots, then polynomials $P^\pm$ also have no real roots; if $R$ has $0 \leq K \leq \frac{n-2}{2}$ real positive roots, then polynomials $P^\pm$ also have $K$ real positive roots. So, the number of positive roots of polynomials $P^\pm$ is even (or zero). Notice that in this case the equator of Poincaré sphere has no equilibria, and so the phase plane contains a "center" for large $r$.

If $> A^2_{\frac{n-2}{2}} > 0$, then one of the polynomials $P^\pm$, say $P^+$, gets one more positive root, $r^+$. We can verify that equilibrium $E^+(r^+, \varphi^+)$ is a saddle. The number of real positive roots of the second polynomial, $P^-$, does not change. So, both polynomials $P^\pm$ have an odd number of equilibria. The equator of the Poincaré sphere in this case contains nodes of alternating stability, and so the phase plane contains "spider-net" structures for large $r$.

Lastly, if $A^2_{\frac{n-2}{2}} < 0$ and $0 < B < -A^2_{\frac{n-2}{2}}$ then $P^-$ gets one more positive root, corresponding to a center, and the number of real positive roots of the polynomial $P^+$ does not change. Polynomials $P^\pm$ both have an even number of positive roots.

Examples of phase portraits of Hamiltonian systems (1.1) with even $n$ are given in Appendix 2 (Figures A2.1 and A2.3).

*3.3. The Hamiltonian as a generalized Lyapunov function*

A derivative of Hamiltonian function (3.3) with respect to system (2.3) is

$$\frac{dH}{dt} = \frac{\partial H}{\partial r}r' + \frac{\partial H}{\partial \varphi}\varphi' = r^2\left(\epsilon_1 + \sum_{k=1}^{s} A_k^1 r^{2k}\right)\frac{\partial H}{\partial \varphi} = B\frac{r^{n-1}}{n}\left(\epsilon_1 + \sum_{k=1}^{s} A_k^1 r^{2k}\right)\sin(n\varphi)$$

Then

$$\left|\frac{dH}{dt}\right| < \epsilon \left|B\frac{r^{n-1}}{n}\right| \text{ if } \left|\left(\epsilon_1 + \sum_{k=1}^{s} A_k^1 r^{2k}\right)\right| < \epsilon,$$

where $\epsilon$ is a small positive number and $k = 1, \ldots, s$.

Thus, for small enough coefficients $\epsilon_1, A_k^1$, the model can be considered as a small perturbation of the Hamiltonian one that keeps its main properties.

## 4. On general model

### 4.1. Finite equilibria

In the general (non-Hamiltonian) case, the $(x,y)$-coordinates of equilibria are defined by the system

$$F_1(x,y) = 0, \quad F_2(x,y) = 0, \tag{4.1}$$

where $F_1(x,y)$, $F_2(x,y)$ are given in (2.4) and the polar coordinates $(r, \varphi)$ of equilibria are defined by the system

$$P(r, \varphi) = 0, \quad Q(r, \varphi) = 0, \tag{4.2}$$

where $P(r, \varphi), Q(r, \varphi)$ are given in (2.3):

$$P(r, \varphi) = r\left(\varepsilon_1 + \sum_{k=1}^{s} A_k^1 r^{2k} + B\, r^{n-2}\cos(n\varphi)\right),$$

$$Q(r, \varphi) = \varepsilon_2 + \sum_{k=1}^{s} A_k^2 r^{2k} - Br^{n-2}\sin(n\varphi).$$

System (2.3) has equilibrium $O$. It can have equilibria $E_k(r_k, \varphi_k)$, where the coordinate $r = r_k$ is the root of the equation

$$F(r) \equiv \left(\varepsilon_1 + \sum_{k=1}^{s} A_k^1 r^{2k}\right)^2 + \left(\varepsilon_2 + \sum_{k=1}^{s} A_k^2 r^{2k}\right)^2 = B^2 r^{2(n-2)}. \tag{4.3}$$

By solving equation (4.3) together with equation

$$F_r(r) = 0, \tag{4.4}$$

one can find *r*-coordinates of multiple equilibrium points of the model and then, using one of equations (4.2), find the corresponding equilibrium value of $\varphi$. Notice now that generally in non-Hamiltonian case, equilibria $E_k(r_k, \varphi_k)$ can be saddles, stable and unstable nodes/spirals. In practical computations, it is more convenient to search peripheral equilibrium points in a neighborhood of quasi-equilibria described in s.2.3. 1

Model (1.1) can have limit cycles (see also 2.3). It was shown in [4,5,8,9] that equilibrium points of general (non- Hamiltonian) system (1.1) as *n*=4 can appeared outside, inside and on a limit cycle. In [4,5,8,9] the sequences of bifurcations of co-dimension 1 that are realized in the system under variation of parameters $A_k^{1,2}$, $B$ have been found.

Notice that we observed similar behaviors for general non- Hamiltonian model (1.1) as *n*=5, *n*=6 (see Fig.A2.7, A2.8).

*4.2. Stability of equilibria*

In a neighborhood of peripheral equilibrium point, the Jacobian matrix $J(r, \varphi) = \begin{pmatrix} P_r & P_\varphi \\ Q_r & Q_\varphi \end{pmatrix}$

of system (2.3) consists of the elements

$P_r = \sum_{k=1}^{s} 2k A_k^1 r^{2k} + (n-2)B\, r^{n-2} \cos(n\varphi), \qquad P_\varphi = -n\, B r^{n-2} \sin(n\varphi),$

$Q_r = \sum_{k=1}^{s} A_k^2\, 2k r^{2k-1} - (n-2)B r^{n-3} \sin(n\varphi), \qquad Q_\varphi = -n\, B\, r^{n-2} \cos(n\varphi).$

In the Hamiltonian case, peripheral points are centers or saddles for which

$$Trace(J(r,\varphi)) = 0$$
$$Det\big(J(r,\varphi)\big) > 0, Det\big(J(r,\varphi)\big) < 0,$$

correspondingly. In the general case (if coefficients $\varepsilon_1, A_k^1$ are small enough ), saddles remain saddles but centers become spirals. It is possible to verify that these spirals are unstable if $O(0,0)$ is unstable. Domains of repelling of all spirals are divided by separatrices (see Figure A2.7, A2.8).

*4.3. Equilibria at infinity (see Figure 3.1a,b)*

**Proposition 8**. *In the equator of Poincaré sphere, equilibrium points of System (2.3) are alternating stable and unstable nodes for odd $n$; the same is true for even $n$ if $B^2 > (A_s^1)^2 + (A_s^2)^2$. When $B^2 < (A_s^1)^2 + (A_s^2)^2$, the system has no equilibria in the equator of the Poincaré sphere.*

Proof of this Proposition is given in Appendix 1.

## 5. Discussion

In this work, we studied complex Equation (1.1) with $n \geq 4$, which describes behaviors of equivariant vector field, i.e. the vector field that is invariant under rotations by the specific angle $2\pi/n$. More precisely, we worked with the Systems (2.3), (2.4) that are equivalent to (1.1); System (2.3) casts Equation (1.1) in polar coordinates, and System (2.4) casts Equation (1.1) in Descartes coordinates. Notice that Systems (2.3), (2.4) are Hamiltonian for specific (explicitly formulated) values of coefficients.

We constructed and described the phase-parameter portraits of Equation (1.1) as "perturbations" of phase-parameter portraits of its Hamiltonian version. The value of *n* (odd or even) as well as coefficients $B$ and $A_s$ serve as the main bifurcation parameters. The structures of portraits essentially depend on these parameters.

We denoted the principle characteristics of considered portraits, which we referred to as "flower rings"; by definition, *"n-flower ring" is a pattern consisting of n centroids and n saddles together with their separatrices; each center is placed inside the "leaf" composed by a separatrix cycle or separatrix loop* (see Figure 1.1 ).

We showed that the number *s* of flower rings in the model (1.1) is $s \leq \frac{n}{2} - 1$ if *n* is even and $s \leq \frac{n-1}{2} - 1$ if *n* is odd. We described also the possible sequences of appearance and disappearance of flower rings in the model under variations of various parameters.

We compared the qualitative structures in the obtained portraits with key features of the ancient ornamental designs that one can see in historical museums of Crete and Athens. On our opinion, some of phase portraits of the model and the ornamental designs have many common characteristics. For this reason, one can consider the Equation (1.1) as a kind of mathematical "blueprint" for these ornaments.

The ornaments are characterized by rings containing a certain number of points connected by spiral-like lines. A remarkable property of these ornaments is the "dynamical indeterminacy", known in the theory of dynamical systems: a small change in the initial conditions of the "starting point", the lines lead to different points, and so one can visit all points of the ornament with a very small shift of a trajectory point. This peculiarity reminded us the phase-parametric portrait of known "weak resonance" complex differential Equation (1.1) proposed by V. Arnold.

We found main patterns of the phase portraits of this equation; some of them are similar to the ornaments, but the portraits may also contain other patterns that appear with parameter variation. Together with "spiral-like" lines and patterns, the ornaments may have patterns composed of closed cycles around centers. Such kinds of patterns usually arise in a Hamiltonian system.

The "spiral-like" ornamental design may reflect the ancient philosophical idea of harmonic unity of the world, where various aspects of world phenomena are interconnected. According to some authors (see, e.g., [11], ch.2), ancient Greeks did not favor the idea of "progress". In contrast, they believed that initially there existed a perfect "golden age", and mankind, in their development down from the Golden Age, are destined to *degenerate* (Hesiod). This degradation may imply broken connections between different parts of the world, causing it to disintegrate into separate parts. The process may be reflected in ornaments by transition from "spiral-like" to "center-like" patterns; we looked for a corresponding exhibit that shows such a transition and have found it in National Archaeological Museum, Athens, see Fig. A3.1.

Mathematically, this process corresponds to the transition from general Equation (1.1) to Hamiltonian (1.1), (3.1), which happens when particular model coefficients vanish. Variation of model parameters results in the change of the phase -parametric portraits, whereby one can "animate the evolution" of ornaments and perhaps on some level reflect the idea of destruction of the "golden age" in the language of mathematics and art.

**Appendix 1**

**Proof of Propositions 3**

***Proposition 3.*** *Let $n \geq 5$ is odd. Then 1) every polynomials $P^+, P^-$ have no more than $\frac{n-1}{2}$ real positive roots; 2) if polynomial $P^+(P^-)$ has K real roots and J of them are positive, then the polynomial $P^-(P^+)$ also has K real roots and K-J of them are positive; 3) every two equilibria $E_j^+, E_{j+1}^+$ corresponding to consequent roots of polynomial $P^+$ are saddle/center or center/saddle; similar statement is valid for the polynomial $P^-$, 4) let $r^*, r_*$ are maximal and minimal roots among all positive roots of polynomials $P^+$ and $P^-$. Then corresponding equilibria $E^*$ and $E_*$ are saddles.*

Coordinates $(r^\pm, \varphi^\pm)$ of peripheral equilibria $E^\pm$ of Hamiltonian system (3H) satisfy the system (3.4). So

$$\varphi_j^\pm = \pm \frac{\pi}{2n} + \frac{2\pi j}{n}, j = 0, \ldots, n-1 \tag{A.1}$$

and $r^\pm$ are the roots of the polynomials $P^\pm$ where

$$P^+(r) = P_0(r) - Br^{n-2}, P^-(r) = P_0(r) + Br^{n-2}, B > 0, P_0(r) = \epsilon_2 + \sum_{k=1}^{s} A_k^2 r^{2k}. \tag{A.2}$$

Even polynomial $P_0(r)$ of the power 2s, where $s = \frac{n-3}{2}$, has $0 \leq m \leq s$ changes of signs of its coefficients, so it has at most $m$ positive roots (or less than $m$ by an even number) due to Descartes theorem. Then one of the polynomials $P^+, P^-$ has at most $m+1$ changes of signs and the second has at most $m$. So, for $m = s$ one of the polynomials $P^+, P^-$ has at most $\frac{n-1}{2}$ positive roots and the second has at most $\frac{n-3}{2}$. Now we show that the number of saddles in the system is at most $\frac{n-1}{2}$, and the number of centers is at most $\frac{n-3}{2}$.

The second assertion of proposition is obvious (evidently, if $r^*$ is a root of $P^+$ then $-r^*$ is a root of $P^-$). The total number of positive roots of both polynomials is an odd number $k \leq n - 2$. Note also that $P^+(0) = P^-(0) = \epsilon_2 \neq 0$, and the equilibrium $O(r = 0, \varphi = 0)$ is a center of the system. The system is Hamiltonian so its equilibria can be only saddles and centersthat *alternates in each ray* $\varphi = \varphi^{\pm} + \frac{2\pi j}{n}, j = 0, ..., n-1$. Let $r_1 < r_2 < \cdots < r_k$ *are positive roots of polynomials* $P^+, P^-$. *Thus the "closest" to O peripheral equilibria* $E_j^{\,1}(r_1, \varphi_j)$ *j=0....n-1 can be only saddles that together with their separatrices compose separatrix cycle in* $(r, \varphi)$- *plane; the cycle contains inside the point O.*

*Next, the largest root* $r_k$ *of polynomials* $P^{\pm}$ *also corresponds to a saddle because the number k is odd; together with their separatrices the saddles* $E_j^{\,k}(r_k, \varphi_j)$ *j=0....n-1 also compose separatrix cycle in* $(r, \varphi)$- *plane.*

*We have proven that the system has at most* $\frac{n-1}{2}$ *separatrix cycles containing n saddles and at most* $\frac{n-3}{2}$ *flower-rings containing n peripheral centers.*

***Proposition 5***. *Let n≥4 be even, and B>0. Then*

1) *polynomials* $P^+, P^-$ *have no more than* $s = \frac{n-2}{2}$ *real positive roots;*
2) *polynomials* $P^{\pm}$ *have the same number of positive roots as polynomial* $P_0$ *if* $A_{\frac{n-2}{2}}^2 > -B > 0$ *and if* $0 < B < -A_{\frac{n-2}{2}}^2$;
3) *if* $B > \left| A_{\frac{n-2}{2}}^2 \right|$, *then at least one of polynomials* $P^-, P^+$ *has a real positive root* $r^*$ *and this* $r^*$*is r-coordinate of a saddle equilibrium of (3.7);*
4) *every two equilibria* $E_j^+, E_{j+1}^+$ *corresponding to subsequent roots of polynomial* $P^+$ *are saddle/center or center/saddle; similar statement is valid for the polynomial* $P^-$.

**Proof of Proposition 5**

In this case $(r, \varphi) = \epsilon_2 + \sum_{k=1}^{s-1} A_k^2 r^{2k} + (A_{(n-2)/2}^2 - B \sin(n\varphi))r^{n-2}$. After substitution $\sin(n\varphi) = \pm 1$ we get even polynomials $P^{\pm}(r, \varphi)$ of the form

$$P^{\pm}(r, \varphi) = \epsilon_2 + \sum_{k=1}^{s-1} A_k^2 r^{2k} + (A_{\frac{n-2}{2}}^2 \mp B)r^{n-2} = P_0(r) \mp Br^{n-2}, \qquad (A.3)$$

$$P_0(r) = \epsilon_2 + \sum_{k=1}^{s-1} A_k^2 r^{2k} + A_{\frac{n-2}{2}}^2 r^{n-2}, \quad B > 0, \quad s = \frac{n}{2} - 1.$$

Each of even polynomial $P^\pm$ has at most $\frac{n-2}{2}$ real positive roots. So, common number positive roots is at most $n - 2$. The equilibrium $O(0,0)$ is the center in the system, so the "closest" to $O$ peripheral equilibria $E_j^1(r_1, \varphi_j), j = 0, \ldots, n - 1$ (if any) are saddles; together with their separartices they compose $n$-separatrix cycle in $(r, \varphi)$-plane. It has been discussed in the proof of Proposition 4 that saddles and centers are alternating in the $(\varphi_j) - rays, j = 0, \ldots, n - 1$.

Let $0 \leq m \leq \frac{n-2}{2}$ is the number of the changes of signs of coefficients of polynomial $P_0(r)$. It is evident that for $A_{\frac{n-2}{2}}^2 > -B > 0$ and for $0 < B < -A_{\frac{n-2}{2}}^2$ both polynomials $P^\pm$ has the same sighs of coefficients for highest exponent as $P_0(r)$. Then coefficients of each polynomial $P^\pm$ have $m$ changes of signs, and so $m$ or less by even number of positive roots; it means that total number of positive roots in both polynomials is even. Then the number of saddles is equal to the number of peripheral centers. Note, that the points $E_j^*(r^* \varphi_j)$, $j$=0….n-1 with the largest root $r^*$ is a center. The union of centroids corresponding to these centers composes $n$-flower-ring.

For $|B| > A_{\frac{n-2}{2}}^2$ one of polynomial $P^-$, $P^+$ "looses" one change of sign in the sequence of coefficients. Then the number of positive roots of this polynomial becomes odd, so total number of positive roots in both polynomials is odd. Then the number of saddles is more by 1 than the number of peripheral centers. Note, that in this case the points $E_j^*(r^* \varphi_j)$, $j$=0….n-1 with the largest root $r^*$ are saddles (see Figure 3.2).

**Proof of Proposition 8**

*Lemma 1.* System (2.4) in equators of Poincare sphere (i.e., in coordinates $(u, w): x = 1/u, y = w/u$ and $(v, u): x = v/u, y = 1/u$) is equivalent to system (2.3) in coordinates $(\rho = \frac{1}{r}, \varphi)$.

Proof of these statement follows from formulas: $\varphi = \arctan(w) = \text{arccot}(v)$, $\rho = \frac{1}{\sqrt{(x^2 + y^2)}}$,

thus $\{\varphi` = \frac{w`}{(1+w^2)}, \rho` = \frac{u`}{\sqrt{(1+w^2)}}\}$ and $\{\varphi` = -\frac{v`}{(1+v^2)}, \rho` = \frac{u`}{\sqrt{(1+v^2)}}\}$.

According to this lemma we consider system (2.3) in $(\rho, \varphi)$ – coordinates. Changing independent variable: $dt \to \rho^{n-2} d\tau$ we get the system:

$$\rho` = -\rho^2(\epsilon_1 \rho^{n-2} + \sum_{k=1}^{s} A_k^1 \rho^{n-2-2k} + B \cos(n\varphi)) \equiv P^i(\rho, \varphi), \qquad (A.4)$$

$$\varphi` = \epsilon_2 \rho^{n-2} + \sum_{k=1}^{s} A_k^2 \rho^{n-2-2k} - B \sin(n\varphi) \equiv Q^i(\rho, \varphi).$$

The equator of Poincare sphere is defined by $= 0$.

For $n$ odd system (A.4) can be written as

$$\rho` = -\rho(B \cos(n\varphi) + o_1(\rho)) \equiv P^{ii}(\rho, \varphi), \qquad (A.5)$$

$$\varphi` = -B \sin(n\varphi) + o_2(\rho)(\rho, \varphi) \equiv Q^{ii}(\rho, \varphi)$$

where $o_1(\rho), o_2(\rho)$ are polynomials such that $o_1(0) = o_2(0) = 0$.

Thus, equilibria in the equator are $I_j(0, \varphi_j)$, where $\varphi_j$ satisfy to the equation $\sin(n\varphi) = 0$, i.e.

$$\varphi_j = \frac{2\pi j}{n}, j = 0, \ldots, n-1.$$

Notice, that $\cos(n\varphi_j) = \pm 1$ for two neighboring $j$. Jacobian of the considering system

$$J(\rho, \varphi) = \begin{pmatrix} P_\rho^{ii} & P_\varphi^{ii} \\ Q_\rho^{ii} & Q_\varphi^{ii} \end{pmatrix} = \begin{pmatrix} -B\cos(n\varphi) & B\rho\sin(n\varphi) \\ 0 & -Bn\cos(n\varphi) \end{pmatrix}.$$

Then $J(0, \varphi_j) = \begin{pmatrix} -B & 0 \\ 0 & -Bn \end{pmatrix}$, $J(0, \varphi_{j+1}) = \begin{pmatrix} B & 0 \\ 0 & Bn \end{pmatrix}$.

Thus, one of these points is a stable node, and another is unstable node. The first statement is proven.

If $n$ is even, then system (A.4) in the Poincare` sphere equators can be written in the form

$$\rho` = -\rho(A_{\frac{n-2}{2}}^1 + B \cos(n\varphi)\rho + o_1(\rho)) \equiv P^{ii}(\rho, \varphi), \qquad (A.6)$$

$$\varphi` = ( A^2_{\frac{n-2}{2}} - B\sin(n\varphi)) + o_2(\rho) \equiv Q^{ii}(\rho, \varphi).$$

Equilibrium values of $\varphi$ in system (A.17) (if they exist) satisfy to the equation

$$A^2_{\frac{n-2}{2}} - B\sin(n\varphi) = 0.$$

This equation has real roots $\varphi_j = \frac{(-1)^j}{n}\arcsin(A^2_{\frac{n-2}{2}}/B) + \frac{\pi j}{n}, j = 0, \ldots, n-1$ if

$\left|A^2_{\frac{n-2}{2}}/B\right| < 1$ and has no roots if $\left|A^2_{\frac{n-2}{2}}/B\right| > 1$.

It is easily to verify that in the second case the equator does not contain equilibria, and in the first case equilibria are alternated stable and unstable nodes. Indeed,

$$J(0, \varphi_j) = \begin{pmatrix} -(A^1_{\frac{n-2}{2}} + B\cos(n\varphi_j)) & 0 \\ 0 & -Bn\cos(n\varphi_j) \end{pmatrix} =$$

$$\begin{pmatrix} -(A^1_{\frac{n-2}{2}} \pm \sqrt{B^2 - \left(A^2_{\frac{n-2}{2}}\right)^2}) & 0 \\ 0 & \mp n\left(A^2_{\frac{n-2}{2}}\right) \end{pmatrix}.$$

*Statements are proven.*

***Corollary.*** Both statements of the Proposition are true also for Hamiltonian systems (4H). Thus, Propositions 4 and 6 are also proven.

**Appendix 2.** *Phase Portraits of several models*

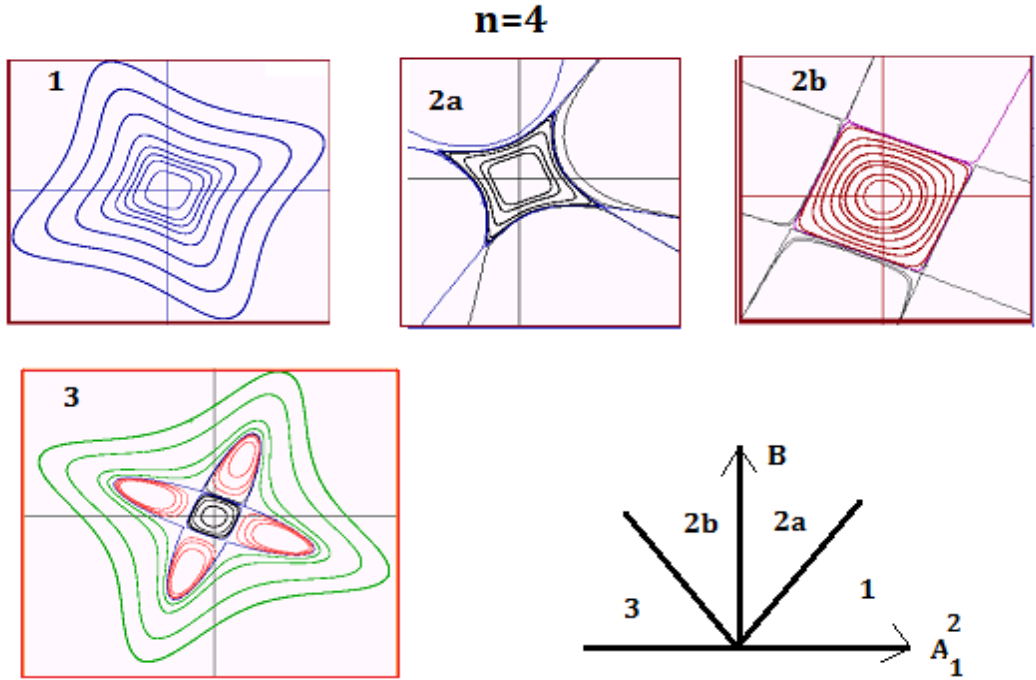

Figure A2.1. Phase-parameter portraits of Hamiltonian model (1) for **n=4**, $\varepsilon_1 = 0, A_1^1 = 0, \varepsilon_2 = 1$.

**1**-"center" ($A_1^2=1, B=.7$), **2a**, **2b**-"center +spider-net" ($A_1^2=1, B=1.2; A_1^2=-1, B=1.2$), **3**- "flower ring+spider-net" ($\varepsilon_2=1, A_1^2=-1, B_2=.7$) The lines $B = A_1^2$, $B = -A_1^2$ are parameter boundaries between domains in the parameters

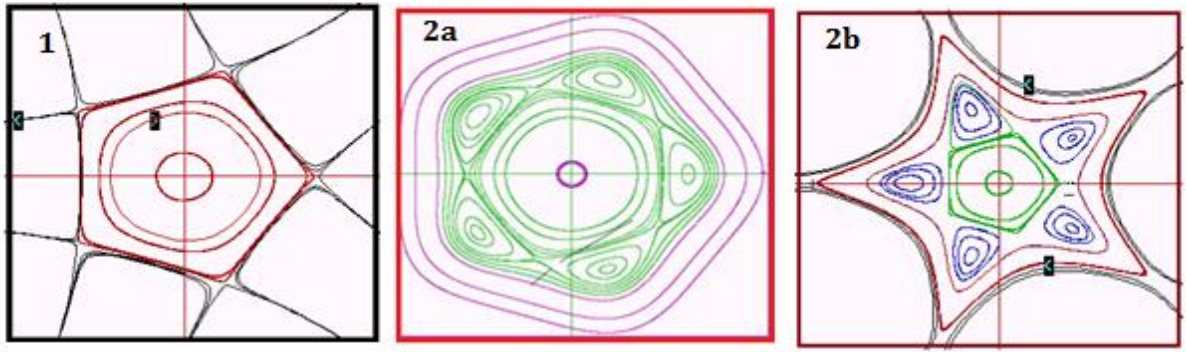

Figure A2.2. Portraits of Hamiltonian model (1.1) for **n=5**, $\varepsilon_1 = 0, A_1^1 = 0$.

**1**-"center+spider-net" ($\varepsilon_2 = 1, A_1^2 = 1, B = 2$), **2a,2b** "flower ring +star+ spider-net"( $\varepsilon_2 = 1, A_1^2 = -1, B = .3$) , ( $\varepsilon_2 = 1, A_1^2 = -1, B = .2$). Phase portraits in Fig.-s **2a** and **2b** are topologically equivalent. The parameter boundary between domains **1** and **2** has equation: $C: 27\epsilon_2 B^2 + 4(A_1^2)^3 = 0$.

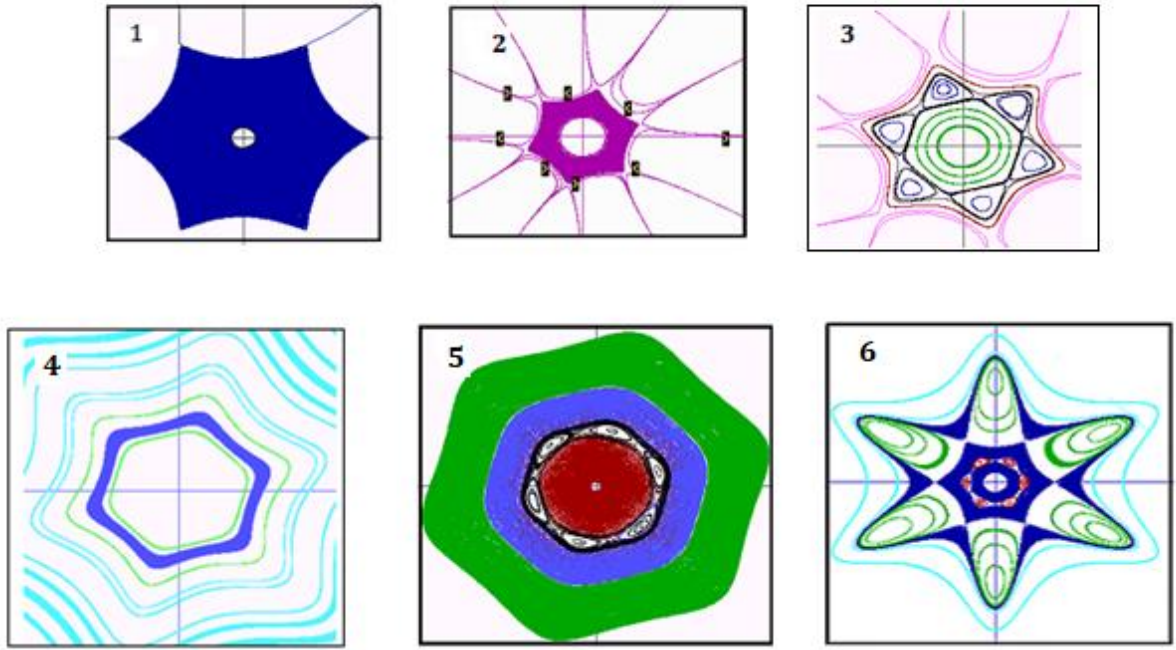

Figure A2.3 Portraits of Hamiltonian model (1) for **n=6,** $\varepsilon_1 = 0, \varepsilon_2 = 1, A_1^1 = 0, A_2^1 = 0.$

**1a, 1b**-"spider-net", $A_1^2 = 1, A_2^2 = 0, B_1 = 0, B_2 = -.5; \varepsilon_1 = \varepsilon_2 = 0, B = .5$ ; **2** –"centers + flower band" , $A_1^2 = -1, A_2^2 = .1, B = .04;$ **3**- "center + star + *two* flower rings", $A_1^2 = -1, A_2^2 = .1, B = .06.$

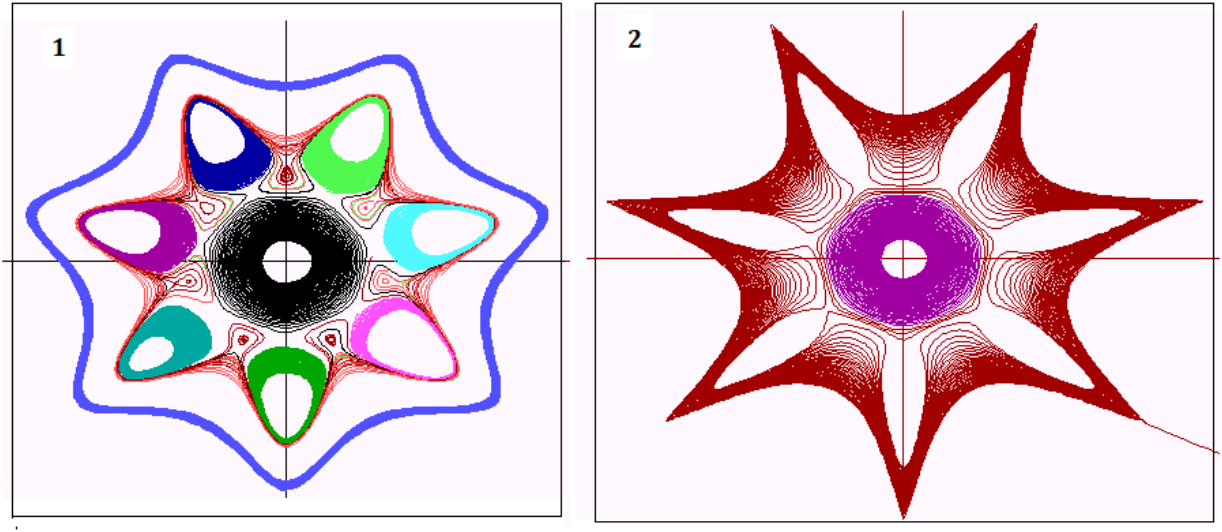

Figure A2.4. Portraits of Hamiltonian model (1) for **n=7**; $\epsilon_1 = A_1^1 = A_2^1 = A_3^1 = 0, \epsilon_2 = -.56$; (**1**) $A_1^2 = 3, \ A_2^2 = -3.5, \ B = -1.6$;(2) $A_1^2 = 3, \ A_2^2 = -3.5, \ B = -1$.

n=9

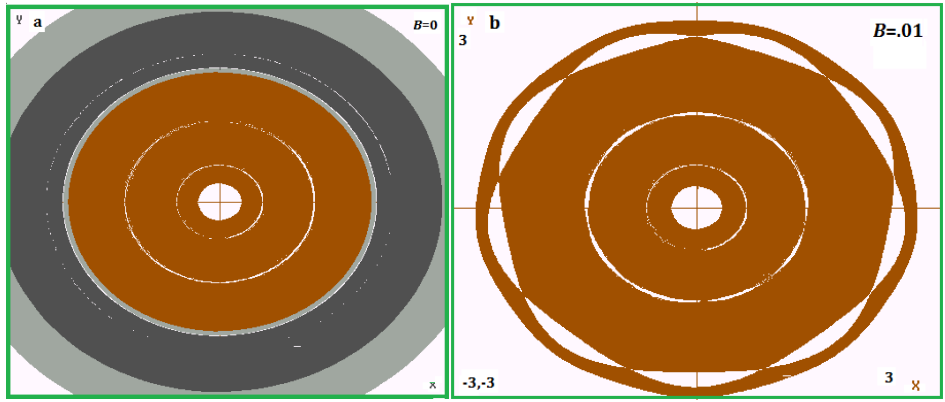

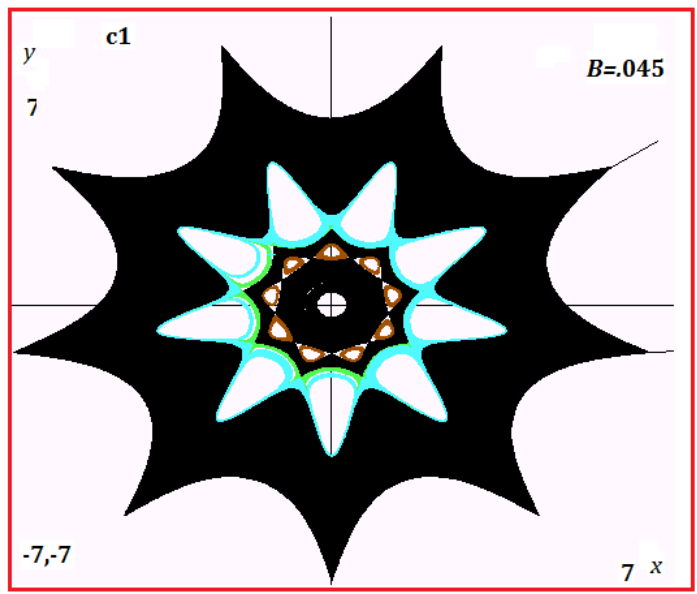

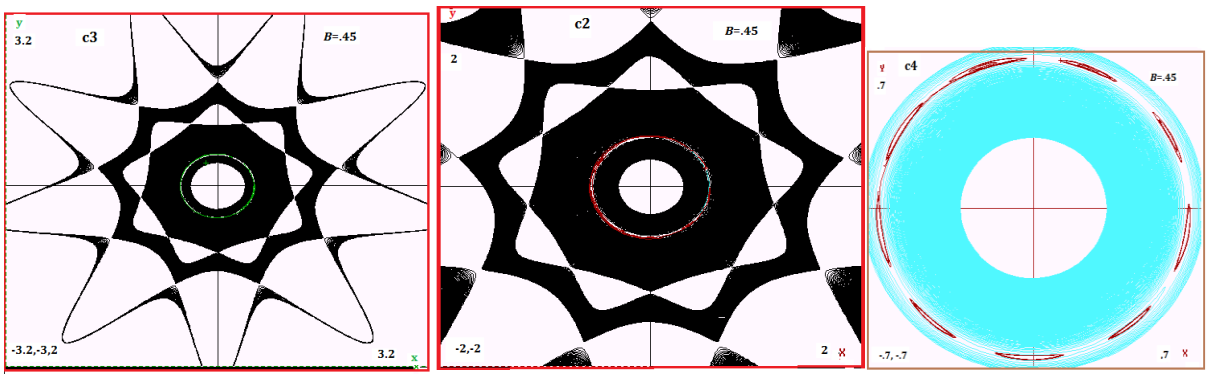

Figure A2.5. Portraits of Hamiltonian model (1) for **n=9**, $\epsilon_1 = A_1^1 = A_2^1 = A_3^1 = 0, \epsilon_2 = 8, A_1^2 = -33, A_2^2 = 23.765, A_3^2 = -3.5.$ **a-**$B$=0, **b-** $B$=.1, **c-** $B$=.45: **windows** in **c1**-$[-7,7]^2$, **c3**-$[-3.2,3.2]^2$, **c2**-$[-2,2]^2$, **c4**-$[-.7,.7]^2$ the picture presents the inside of center-ring that does not seen clearly in Fig.-s **c3**(the green cycle), **c2** (the red cycle)

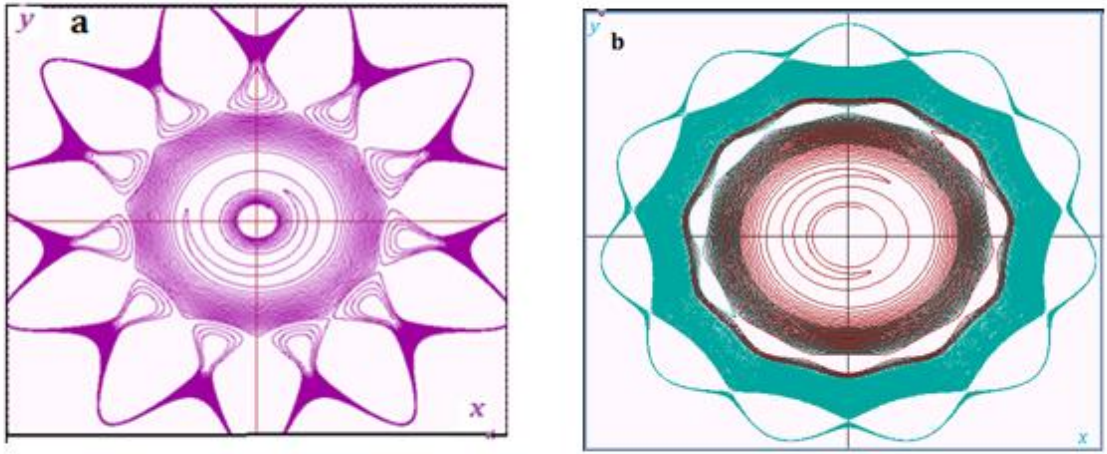

Figure A2.6. Portraits of Hamiltonian model (1) for **n=11**. $\epsilon_1 = A_1^1 = A_2^1 = A_3^1 = A_4^1 = 0, \epsilon_2 = 14.4$, $A_1^2 = -55.6$, $A_2^2 = 54.6$, $A_3^2 = 14.4$, $A_4^2 = 1$; two flower rings. **a** – the portrait for $B = .05$, **b** – the portrait for $B = -.01$

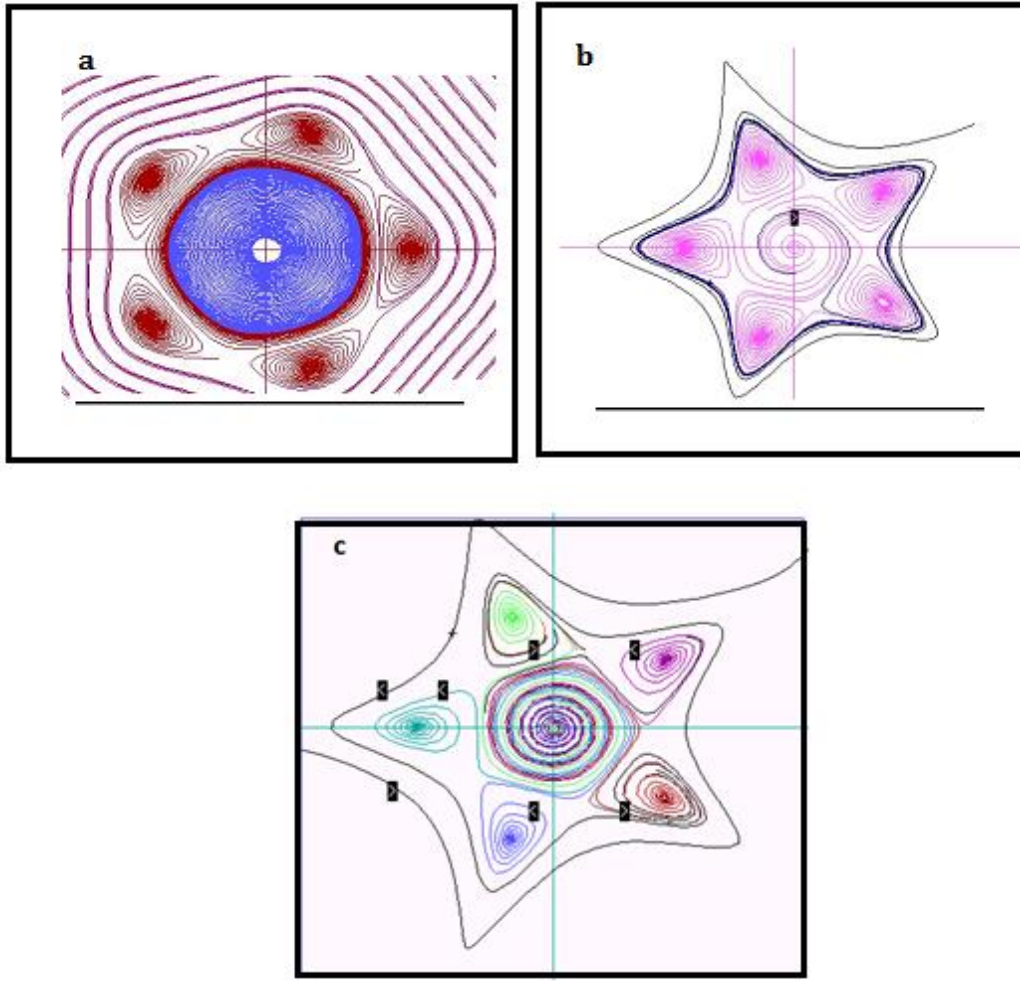

Figure A2.7. Portraits of non- Hamiltonian model (1) for **n=5**, $\epsilon_1 = .005$.

" Flower ring" outside (**a**) and inside (**b**) of limit cycle, (**c**) without limit cycle,

**a**: $\epsilon_2 = 1, A_1^1 = -.01, A_1^2 = -1, B = .1$; **b:** $\epsilon_2 = -.1, A_1^1 = .045, A_1^2 = 1, B = 1$;

**c :** $\epsilon_2 = -.1, A_1^1 = -.045, A_1^2 = 1, B = 1$.

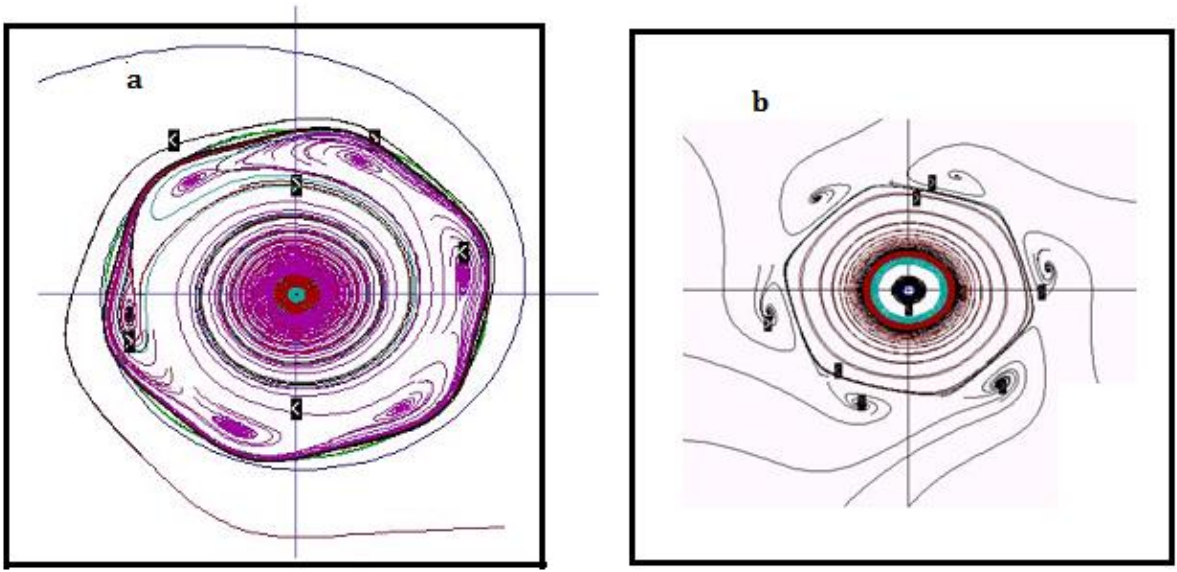

Figure A2.8. Portraits of non- Hamiltonian model (1) for **n=6**, $\epsilon_1 = -.001$.

" Flower ring" peripheral equilibria are placed outside (**a**) and inside (**b**) of unstable limit cycle,

**a**: **a:** $\epsilon_2 = -.1,\ A_1^1 = 1.3,\ A_1^2 = .1,\ A_2^2 = -.1, B = .05$

**b:** , $\epsilon_2 = .1,\ A_1^1 = 1,\ A_1^2 = -.1,\ A_2^2 = -.1, B = .05$

**Appendix 3.** *Designs presented in museums of Crete and Athens*

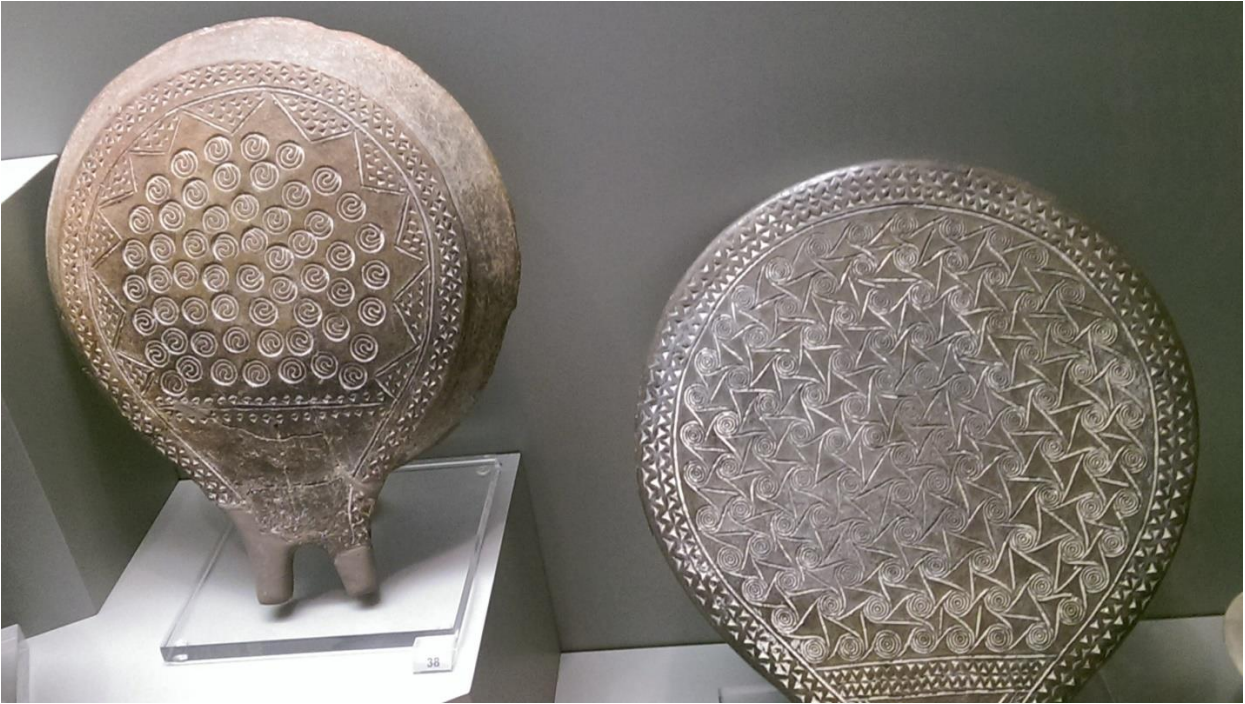

Figure A3.1

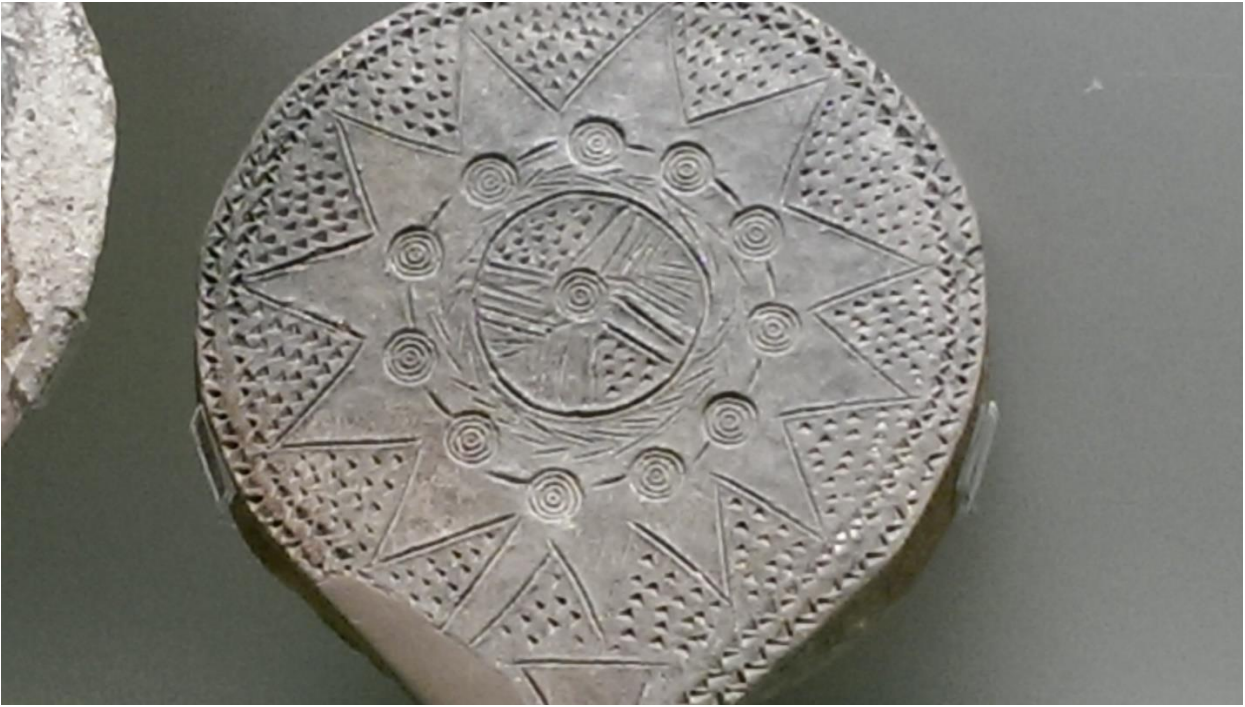

Figure A3.2

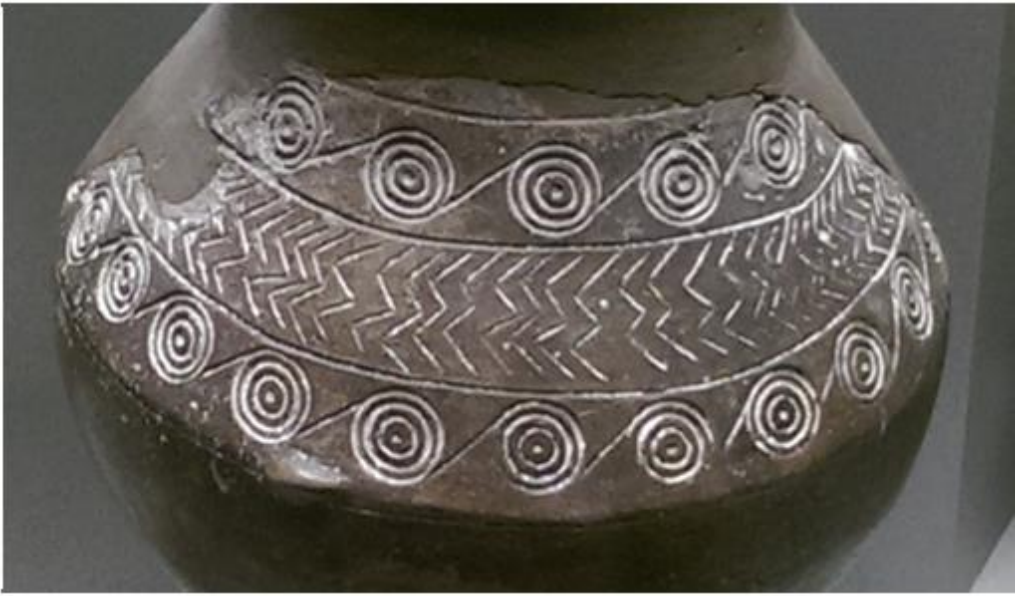

Figure A3.3

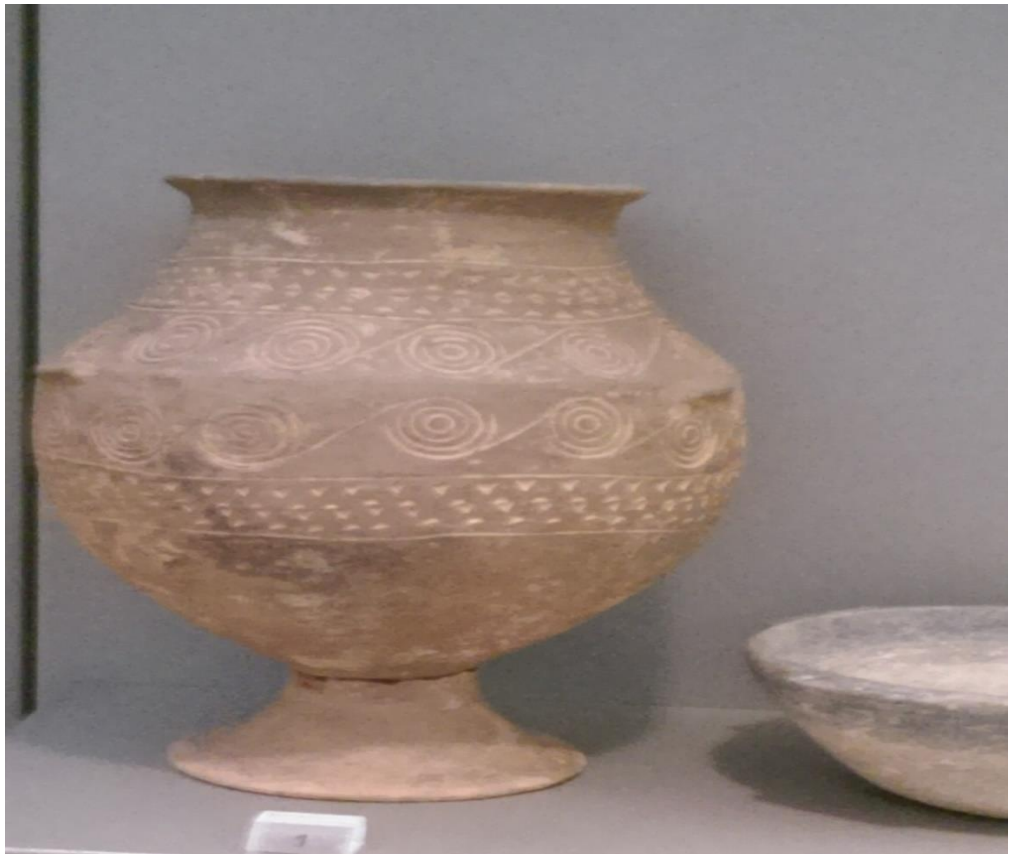

Figure A3.4

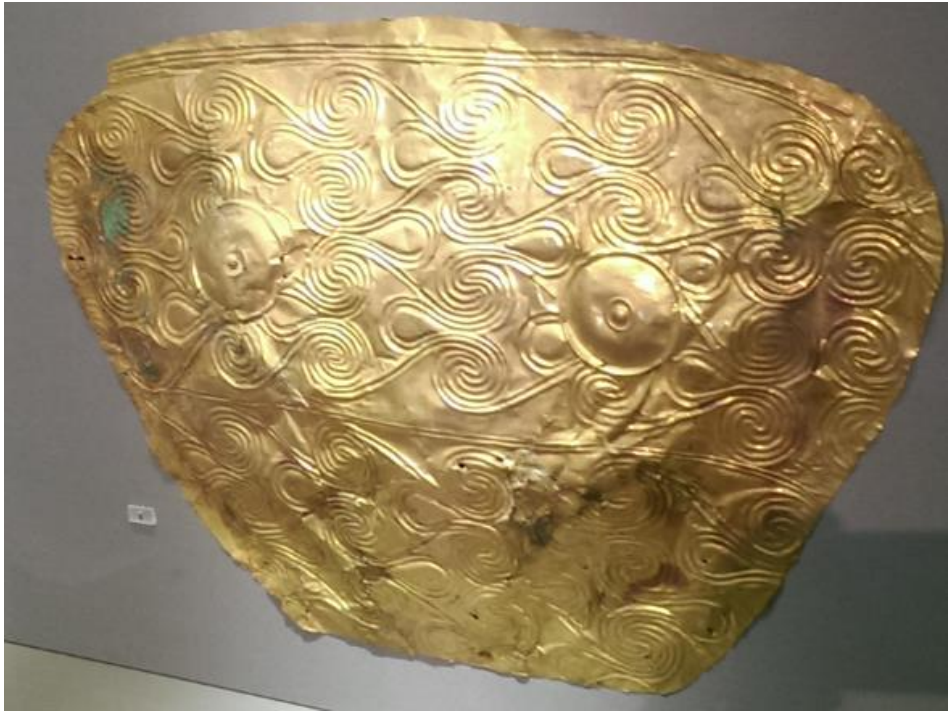

Figure A3.5

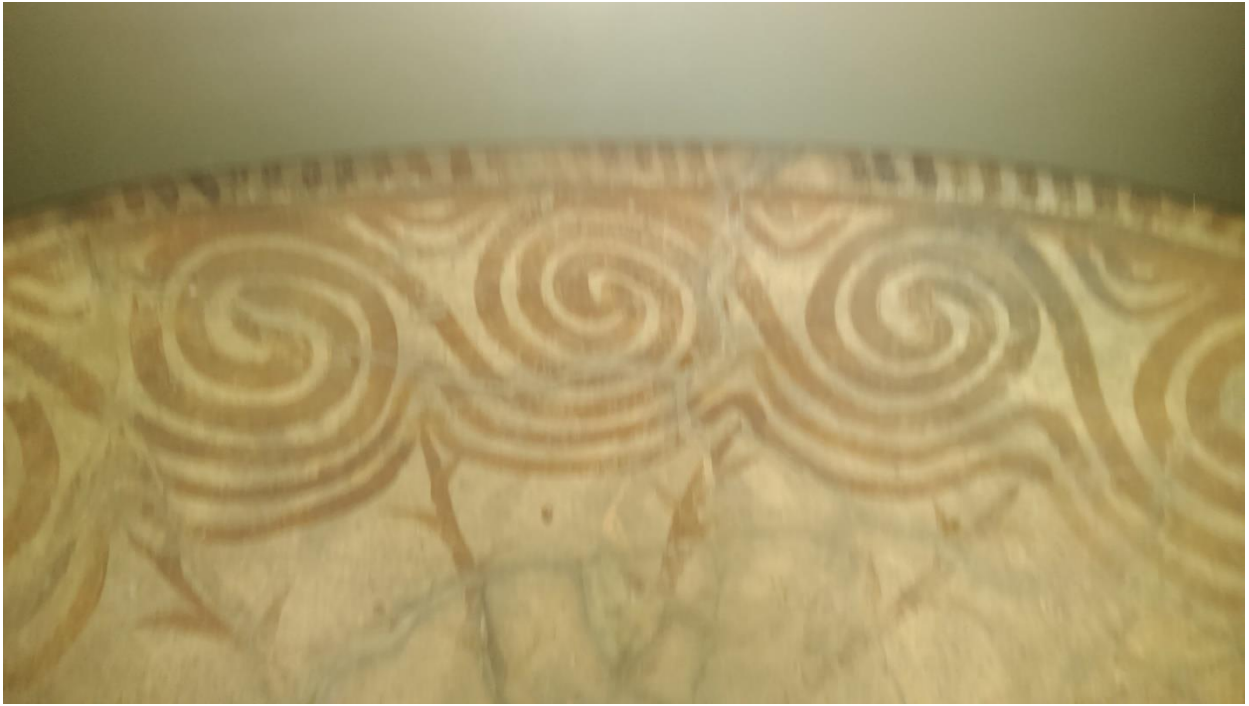

Figure A3.6